\documentclass[reqno]{amsart}
\usepackage{amsmath, amssymb, amsthm, geometry, enumerate, graphicx, amsfonts, hyperref, mathrsfs}
\hypersetup{colorlinks=true,linkcolor=red, anchorcolor=blue, citecolor=blue, urlcolor=red, filecolor=magenta, pdftoolbar=true}
\theoremstyle{plain}
\newtheorem{thm}{\it Theorem}[section]
\newtheorem{prop}[thm]{\it Proposition}

\newtheorem{lem}[thm]{\it Lemma}
\theoremstyle{remark}
\newtheorem{defn}[thm]{Def{}inition}
\newtheorem{rem}[thm]{Remark}
\newtheorem{exa}[thm]{Example}
\numberwithin{equation}{section}

\usepackage{enumerate}
\usepackage{relsize}
\allowdisplaybreaks
\begin{document}

\title [Frames for the Weyl--Heisenberg group and the extended affine group]{Nonstationary frames of translates and frames  from  the Weyl--Heisenberg group and the extended affine group }
\author[Divya Jindal]{Divya Jindal}
\address{{\bf{Divya Jindal}}, Department of Mathematics,
University of Delhi, Delhi-110007, India.}
\email{divyajindal193@gmail.com}

\author[Lalit   Kumar Vashisht]{Lalit  Kumar  Vashisht$^*$}
\address{{\bf{Lalit  Kumar  Vashisht}}, Department of Mathematics,
University of Delhi, Delhi-110007, India.}
\email{lalitkvashisht@gmail.com}

\begin{abstract}
In this work, we analyze Gabor  frames for  the Weyl--Heisenberg group and wavelet frames for the extended affine group. Firstly, we give necessary and sufficient conditions for the existence of nonstationary frames of translates. Using these conditions, we give the existence of Gabor frames from the Weyl--Heisenberg group and wavelet frames for the extended affine group. We present a representation of functions in the closure of the linear span of a Gabor frame sequence in terms of the Fourier transform of window functions.  We show that the canonical dual of frames of translates has the same structure.  An approximation of inverse of the frame operator of nonstationary frames of translates is presented. It is shown  that a nonstationary frame of translates is a Riesz basis if it is linearly independent and satisfies  approximation  of the inverse  frame operator. Finally, we give equivalent conditions for a nonstationary sequence of translates to be linearly independent.
\end{abstract}

\renewcommand{\thefootnote}{}

\footnote{
The research of Divya Jindal is supported by the Council of Scientific $\&$ Industrial Research (CSIR), India.  Grant No.: 09/045(1680)/2019-EMR-I. Lalit Kumar Vashisht is  supported by the Faculty Research Programme Grant-IoE, University of Delhi \ (Grant No.: Ref. No./IoE/2021/12/FRP).\\
$^*$Corresponding author: Lalit Kumar Vashisht.}

\maketitle

\baselineskip15pt
\section{Introduction}
Frames for separable Hilbert separable spaces were introduced by  Duffin and Schaeffer in \cite{DS}; originally based on work by Gabor \cite{Gabor} in decomposition of signals in terms of elementary functions. The fundamental concept of frames was revived by  Daubechies, Grossmann, and Meyer \cite{DGM}, who showed its importance for data processing.  Frames were also reviewed by Young in \cite{Young}. Since then, the theory of frames began to be studied more widely.   Wavelet analysis and Gabor analysis are part of this development, they are an integral part of  time-frequency analysis. Wavelet analysis corresponds to unitary irreducible representations of the affine group, while  Gabor analysis corresponds to unitary irreducible representations of the Weyl-Heisenberg group. Whereas Gabor analysis yields a time-frequency representation of signals, wavelet analysis provides a time-scale representation. There is huge literature on  Gabor analysis and   wavelet analysis and its applications in both pure mathematics and engineering science, see \cite{ole, DabI,   KG, Bhan, H89, Heil20, Hern1, Jan, PalleyJ, AK, lopez, Meyer, Hari, Woj} and  many references therein.

In 1998, Dai and  Liang proved that if $\psi_1$ and $\psi_2$ are two affine orthonormal  MRA wavelets then there exists a continuous map $A:  [0,1] \rightarrow L^2(\mathbb{R}, dx)$ such that $A(0)=\psi_1$, $A(1)=\psi_2$, and $A(t)$ is an affine orthonormal MRA wavelet for all $t \in [0,1]$. The proof appeared in \cite[Theorem 4]{Wu}, and may also be found in the review article \cite{WW} by Weiss and Wilson. Years later, Dahlke,  Fornasier,  Rauhut,  Steidl and  Teschke in \ \cite{DFR} and Torresani in \cite{BT1,BT2} established a connection  between  Gabor analysis and wavelet analysis by exhibiting both the affine group and the  Weyl--Heisenberg group  as subgroups of the four dimensional affine Weyl--Heisenberg group.  They also contracted the corresponding coherent state families, resolutions of the identity, and tight frames.  One of the key insights in this paper is the realization that the affine group and the Weyl--Heisenberg group have different dimensions. They overcame this problem by introducing the extended affine group and an additional parameter, $\varepsilon$. The frames thus obtained vary continuously with $\varepsilon$, $\varepsilon=0$ yields a Gabor frame and $\varepsilon=1$ yields a wavelet frame.

Inspired by the above works we study the characterization of Gabor   frames in $L^2(\mathbb{R}, dx)$ with compact support and several generators (and in particular, Parseval frames). We apply our results in the context of group representations, and in particular to those systems associated with the Weyl--Heisenberg group and extended affine group, using the unitary irreducible representation. Whereas authors of  \cite{DGM} and \cite{SB} obtained sufficient conditions, we find explicit necessary and sufficient conditions for generators having a compact support. Recently, the authors in  \cite{JVas} proved sufficient conditions for the existence of Gabor frames  and wavelet frames with several generators for the Weyl--Heisenberg group and the extended affine group, respectively. They also proved a Paley-Wiener type perturbation results for frames with several generators.  Extension of frame conditions from $L^2(\mathbb{R}, dx)$ to its corresponding matrix-valued Lebesgue spaces and vice versa were studied  in \cite{JVIII}.   Notable contribution in this paper includes necessary and sufficient conditions for nonstationary frames of translates, Gabor frames and wavelet frames with several generators associated with the Weyl--Heisenberg group and the extended affine group.  Using results of  Casazza and  Christensen \cite{Casa} and    Christensen and Hasannasab \cite{olehasan}, we give  approximation of inverse frame operator and linear independence of nonstationary frames of the space $L^2(\mathbb{R}, dx)$.

The work in this paper is structured as follows. In Section \ref{SecI}, we  briefly review  Hilbert frames, Riesz basis, Gabor systems and wavelet systems. The main results start from Section \ref{Sec2}. Theorem \ref{w1} gives necessary and sufficient conditions  for the existence of Gabor frames associated with the Weyl--Heisenberg group. Its proof is based on the existence of nonstationary frames of translates or generalized shift-invariant frames for $L^2(\mathbb{R},dx)$, see Theorem \ref{t1}. In Theorem \ref{3.77}, we give an interplay between modulation and translation parameters in Gabor frames associated with the Weyl--Heisenberg group. Theorem \ref{thm3.10x} gives a representation of a function in the closure of linear span of a Gabor frame sequence in terms of Fourier transform of window functions.   In Section \ref{Sec4}, we give necessary and sufficient conditions for  wavelet frames  associated with the extended affine group, see Theorem \ref{th4.8e}.  In Section \ref{Sec5}, we discuss the structure of the canonical dual of nonstationary frames of translates for $L^2(\mathbb{R},dx)$. In Theorem \ref{th5.1}, we show that the canonical dual of  nonstationary frames of translates has the same structure. By using a technique given in \cite{Casa}, we give an approximation of inverse of the frame operator of nonstationary frames of translates in Theorem \ref{r1} of  Section \ref{Sec6}. Sufficient condition for nonstationary Riesz bases of the space $L^2(\mathbb{R},dx)$ are given in Theorem \ref{r3}. In Section \ref{Sec7}, we give equivalent conditions for linear independence of a nonstationary sequence of translates in $L^2(\mathbb{R},dx)$.

\section{Preliminaries}\label{SecI}
As usual $\mathbb{N}$, $\mathbb{Z}$, $\mathbb{R}^{+}$, $\mathbb{R}^{*}$, $\mathbb{R}$ denote the set of positive natural numbers, integers, positive real numbers, non-zero real numbers and real numbers, reepectively. $\mathbb{I}$ denotes a countable (finite or infinite) index set. The \emph{support} of a function $f$ defined on $\mathbb{R}$  is closure of the set $ \{x \in \mathbb{R}: f(x) \ne 0\}$.  Throughout this paper, $\mathcal{H} \ne \{0\}$ denotes a  separable Hilbert space with  inner product $\langle \cdot , \cdot \rangle$ and standard norm on $\mathcal{H}$ is given by
$\|f\| = \sqrt{\langle f , f\rangle}$, $f \in \mathcal{H}$.
\newpage
\subsection {Hilbert frames and Riesz Bases:}  A collection of vectors $\mathcal{F}:=\{ \varphi_k \}_{k \in \mathbb{I}}$ in  $\mathcal{H}$ is  called a \emph{frame} (or \emph{Hilbert frame}) of  $\mathcal{H}$ if for  some $\alpha_o$, $\beta_o \in (0, \infty)$  the following inequality holds
\begin{align}\label{fineq}
\alpha_o \|f\|^2\leq  \sum\limits_{k \in \mathbb{I}} |\langle f , \varphi_k\rangle|^2 \leq \beta_o \|f\|^2
\end{align}
for all $f \in \mathcal{H}$. Inequality  \eqref{fineq} is called the \emph{frame inequality}. The positive scalars $\alpha_o$ and $\beta_o$, obviously not unique,  are known as \emph{lower frame bound} and \emph{upper frame bound} of $\mathcal{F}$.  The frame $\mathcal{F}$ is \emph{ tight} if $\alpha_o =\beta_o$, it is \emph{Parseval},  if  $\alpha_o =\beta_o =1$. If $\mathcal{F}$   fulfills the upper inequality  in \eqref{fineq},   then we say that  $\mathcal{F}$ is  a \emph{Bessel sequence} with \emph{Bessel bound} $\beta_o$. If $\mathcal{F}$ is a Bessel sequence, then the map $S: \mathcal{H} \rightarrow \mathcal{H}$ given by $Sf = \sum\limits_{k \in \mathbb{I}} \langle f, \varphi_k \rangle \varphi_k$ is called the \emph{frame operator} of the frame $\mathcal{F}$. The frame operator $S$ is bounded and linear. It is invertible on $\mathcal{H}$ if $\mathcal{F}$ is a frame for $\mathcal{H}$. In that case $\mathcal{F}$ gives the following reconstruction formula,
  $f = SS^{-1}f =\sum\limits_{k \in \mathbb{I}} \langle f , S^{-1}\varphi_k \rangle \varphi_k$,  $f \in \mathcal{H}$. The scalars $\langle f , S^{-1}\varphi_k \rangle$ are known as \emph{frame coefficients}. The supremum over all  lower frame bounds is called the \emph{optimal lower frame bound},  and the infimum over all upper frame bounds is called the \emph{optimal upper frame bound}. The optimal bounds in terms of the frame operator are given in the following result.
\begin{prop}\label{o1}\cite[p. 121]{ole}
	The lower optimal frame bound $\alpha_{\text{Opt}}$ and the upper optimal frame bound $\beta_{\text{Opt}}$ of a frame $\mathcal{F}$ with frame operator $S$ are given by $\alpha_{\text{Opt}}= \|S^{-1}\|^{-1}$, $\beta_{\text{Opt}}= \|S\|$.
\end{prop}

\begin{defn}\cite[Definition 7.3.1]{ole}
A sequence $\{\varphi_k \}_{k \in \mathbb{I}}$ in $\mathcal{H}$ is a \emph{Riesz frame} of $\mathcal{H}$ if every subsequence of $\{\varphi_k\}_{k \in \mathbb{I}}$ is a frame for its closed linear span, with the same frame bounds $\alpha$, $\beta$.
\end{defn}

\begin{defn}\cite[Definition 3.6.1]{ole}
Let  $\Xi$ be a bounded, linear and bijective operator acting on  $\mathcal{H}$   and $\{\chi_k\}_{k \in \mathbb{I}}$ be an orthonormal basis of $\mathcal{H}$.  A sequence of the form $\{\Xi \chi_k\}_{k \in \mathbb{I}}$ is called a \emph{Riesz basis} for $\mathcal{H}$.
\end{defn}
\begin{thm}\cite[Theorem 7.1.1]{ole}\label{Riesb}
A frame  $\mathcal{F} = \{\varphi_k\}_{k\in \mathbb{I}}$ of $\mathcal{H}$ is a Riesz basis for $\mathcal{H}$ if and only if   $\mathcal{F}$ has a biorthogonal sequence $\{\psi_k\}_{k\in \mathbb{I}}$, that is,
\begin{align*}
\langle \varphi_i, \psi_j \rangle= \delta_{ij}=\begin{cases}
1, \quad \text{if}\ i=j;\\
0, \quad \text{elsewhere}.
\end{cases}
\end{align*}
\end{thm}
\begin{rem}\label{remRies}
Every Riesz basis is a basis and hence linearly independent. For other types of independence of frames we refer to Chapter 7 of \cite{ole}.
\end{rem}
The frame conditions given in \eqref{fineq} is  a powerful tool in the study of operator theory \cite{Aldr, JytIII, Meyer}, iterated function systems \cite{D, VD3}, quantum physics \cite{AAG,  JytI, JytII}, distributed signal processing \cite{DVII, DVIII}. Among the many available texts on frames,  Christensen  \cite{ole},  Heil \cite{Heil20} are excellent for basic theory of frames and Han \cite{Bhan} for applications of frames in many directions.

\subsection{Gabor System and Wavelet System}
As is standard, $L^2(\mathbb{R},dx): = \Big\{f: \int_{\mathbb{R}} |f|^2 dx < \infty \Big\}$, the space of  square integrable (in the sense of Lebesgue) functions  over $\mathbb{R}$. $L^2(\mathbb{R},dx)$ is a Hilbert space with respect to the standard inner product
$\langle f, g \rangle = \int_{\mathbb{R}} f(x) \overline{g(x)} dx$.
  The Fourier transform of a function $f$, denoted by $\widehat{f}$, defined as
\begin{align*}
\widehat{f}(\gamma) =  \int_{\mathbb{R}} f(x) e^{- 2\pi i x \cdot \gamma} dx, \ \ \gamma \in \mathbb{R}.
\end{align*}
Mathematically, a Gabor system and wavelet system are defined by using the following three classes of operators which act unitarily on $L^2(\mathbb{R},dx)$. For $a$, $b\in \mathbb{R}$ and  $c \in \mathbb{R}^*$ define translation operator, modulation operator and dilation operator, respectively,  on $L^2(\mathbb{R},dx)$ by
\begin{align*}
T_a f(x) &\mapsto f(x-a), \ f \in L^2(\mathbb{R},dx)  \quad \text{(Translation by a);}\\
E_b f(x) &\mapsto e^{2\pi ibx}f(x), \ f \in L^2(\mathbb{R},dx)  \quad \text{(Modulation by b);}\\
D_c  f(x) &\mapsto \frac{1}{\sqrt{c}}f\Big(\frac{x}{c}\Big), \ f \in L^2(\mathbb{R},dx)   \quad \text{(Dilation by c).}
\end{align*}
Let $\phi \in L^2(\mathbb{R},dx)$ be a non-zero function. A collection   of functions of the form $\mathcal{G}(a, b, \phi): = \{E_{mb}T_{na}\phi\}_{m, n \in \mathbb{Z}}$  $= \{e^{2 \pi i m b x} \phi(x -na)\}_{m, n \in \mathbb{Z}}$ is called the \emph{Gabor system}; and the collection $\mathcal{W}(c, b, \phi): \{T_{k b c^j}D_{c^j}\phi\}_{j, k \in \mathbb{Z}} = \{c^{\frac{-j}{2}} \phi(c^{-j} x- k b)\}_{j, k \in \mathbb{Z}}$ is called the \emph{wavelet system}.  A frame of the form  $\mathcal{G}(a, b, \phi)$ and $\mathcal{W}(a, b, \phi)$ for $L^2(\mathbb{R},dx)$ is  called  the \emph{Gabor frame} and the  \emph{wavelet frame}, respectively.  Christensen \cite{ole} and Heil \cite{Heil20} are good texts for fundamental properties of Gabor frames and wavelet frames.

The following  fundamental properties of above operators will be used throughout the paper.
\begin{lem}\cite[p. 65]{ole}
For any  $0\neq a \in L^2(\mathbb{R},dx)$ and $f \in L^2(\mathbb{R})$, we have
\begin{align*}
\widehat{T_af}=E_{-a}\widehat{f}, \quad \widehat{E_a f}=T_a \widehat{f}, \quad \widehat{D_a f}=D_{\frac{1}{a}}\widehat{f}.
\end{align*}
\end{lem}
\begin{lem}\cite[p. 65]{ole}\label{p1}
	For $f \in L^2(\mathbb{R}), a,b \in \mathbb{R}, c \in \mathbb{R}^*$  and $x \in \mathbb{R}$, the following commutator relations hold.
\begin{align*}
T_aE_bf(x)=e^{-2\pi i ba}E_bT_af(x), \quad T_aD_cf(x)=D_cT_{a/c}f(x), \quad D_cE_bf(x)=E_{b/c}D_cf(x).
\end{align*}
\end{lem}
The following result is direct consequence of Parseval's equation and can be found in any standard text on analysis, for instance \cite{Heil20, Young}.
\begin{thm}\label{thmpp1}
 Assume that for $\mu>0$ and $f \in L^2([0,\mu],dx) \subset L^2(\mathbb{R},dx)$, we have
\begin{align*}
\mu \|f\|_{L^2(\mathbb{R},dx)}^2=\sum_{k \in \mathbb{Z}}|c_k|^2, \quad \text{where} \ c_k=\int_{\mathbb{R}}f(x)e^{-2 \pi i\frac{k}{\mu}x} dx.
\end{align*}
\end{thm}
\subsection{Nonstationary Frames of Translates}
Motivated by the concept of nonstationary wavelet system in  \cite[ Chapter  8]{AK}, we consider  nonstationary frames of translates.  Let $X$ be a countably infinite index set 
and for each $j \in X$, let $\varphi_j \in L^2(\mathbb{R})$. In this work, we consider the following nonstationary frames of translates:
\begin{defn}
Let  $a \in \mathbb{R}$ and $\varphi_j \in L^2(\mathbb{R})$, $j \in \mathbb{Z}$.  A frame of  the form $\{T_{k a^{j}} \varphi_j \big\}_{j, k \in \mathbb{Z}}$  for $L^2(\mathbb{R}, dx)$ is called the  \emph{nonstationary frame of translates} of $L^2(\mathbb{R}, dx)$.
\end{defn}
 In case of nonstationary frames of  translates, we consider infinitely many window functions $\varphi_j$, whereas in  stationary frames of  translates, we consider finitely many window functions. For basic results on stationary Bessel sequences and stationary frames of translates (with one window function), we refer to \cite[Chapter 7]{ole}. 
Nonstationary frames with different structure studied by many authors, we refer to  \cite{AK} for nonstationary  wavelets and related applications. Most recently, Jindal, Jyoti and Vashisht \cite{JVIII} studied matrix-valued nonstationary frames in the matrix-valued signal space $L^2(\mathbb{R}, \mathbb{C}^{n \times n})$.  

\section{Frames  for the  Weyl--Heisenberg Group}\label{Sec2}
Let $\mathcal{G}$ be a locally compact group with left Haar measure $\mu$. A unitary representation of $\mathcal{G}$ is a homomorphism $\pi$ from $\mathcal{G}$ into the group $\mathcal{U} (\mathcal{H})$ of unitary operators on $\mathcal{H}$ that is continuous with respect to the strong operator topology, that is, a map $\zeta: \mathcal{G}\to \mathcal{U}(\mathcal{H})$ that satisfies $\zeta(xy)=\zeta(x)\zeta(y)$ and $\zeta(x^{-1})=\zeta(x)^{-1}=\zeta(x)^*$, and for which $x\to \pi(x)u$ is continuous from $\mathcal{G}$ to $\mathcal{H}$ for any $u\in \mathcal{H}$.  Also, $\zeta$ is said to be irreducible if $\pi$ admits only trivial (means, $=\{0\}$ or $\mathcal{H}$) invariant subspaces.

The Weyl--Heisenberg group $\mathcal{W}$ is the outer semidirect product $\mathbb{R}^2 \rtimes_\Phi \mathbb{R}$, where $\Phi:\mathbb{R} \rightarrow Aut(\mathbb{R}^2)$ given by $\Phi_{\zeta}\bigg(\begin{pmatrix}
a_1 \\
a_2
\end{pmatrix}\bigg) = \begin{pmatrix}
a_1+\zeta a_2 \\
a_2
\end{pmatrix}$ for every $\zeta \in \mathbb{R}$ and $\begin{pmatrix}
a_1 \\a_2
\end{pmatrix} \in \mathbb{R}^2$. The product of $(s,\vec{u})$ and $(t, \vec{v})$ in $\mathcal{W}$ is given by
\begin{align*}
(s,\vec{u})(t, \vec{v})=(s+t, \Phi_{\zeta}(\vec{u}+\vec{v})).
\end{align*}
For every $A\in \mathbb{R}^*, B \in \mathbb{R}$, the unitary irreducible representation of $\mathcal{W}$ is given by
\begin{align*}
\Theta^{A,B} \colon \mathcal{W} &\rightarrow L^2 (\mathbb{R},dx)\\
(\Theta^{A,B} (c,(y_1,y_2))f)(x) &= e^{2\pi i[A(y_1 +xy_2) + By_2]}f(c+x).
\end{align*}
All unitary irreducible representations of this group $\mathcal{W}$ are unitarily equivalent, see \cite{Mackey, Neumann, Neumann1, LT} for technical details.
We refer to the \cite{JVas, SB} for technical details about applications of the  Weyl--Heisenberg group in the construction of Gabor frames.

For $j\in \{1,2,\dots,N\}$, where $N$ is some strictly positive number, let $p_0^{(j)}$ and $q_0^{(j)}$ be real numbers such that $|p_0^{(j)}q_0^{(j)}|<1$, and consider the discrete subset of $\mathcal{W}$ given by
\begin{align*}
\mathcal{W}_{q_0^{(j)},p_0^{(j)}}=\bigg\{\bigg(nq_0^{(j)},\bigg(\frac{nlq_0^{(j)}p_0^{(j)}}{2},lp_0^{(j)}\bigg)\bigg)|n,l \in \mathbb{Z}\bigg\}.
\end{align*}
For each  $j\in \{1,2,\dots,N\}$, let  $\Phi_j$ be a non-zero function in  $L^2(\mathbb{R},dx)$ (also known as ``window function''). For $n,l \in \mathbb{Z}$, define
\begin{align*}
\Phi_{(n,l,j)}^{A,B}=\Theta^{A,B}\bigg(nq_0^{(j)},\bigg(\frac{nlq_0^{(j)}p_0^{(j)}}{2},lp_0^{(j)}\bigg)\bigg)\Phi_j,
\end{align*}
Then
\begin{align*}
\Big\{\Phi_{(n,l,j)}^{A,B}(x)\Big\}_{\underset{j\in \{1,2,\dots,N\}}{n,l\in \mathbb{Z}}} &= \Big\{e^{2\pi i[(1/2)Anlp_0^{(j)}q_0^{(j)}+Blp_0^{(j)}]} e^{2\pi i Alp_0^{(j)}x} \Phi_j(x+nq_0^{(j)})\Big\}_{\underset{j\in \{1,2,\dots,N\}}{n,l\in \mathbb{Z}}}\\
&= \Big\{e^{2\pi i[(1/2)Anlp_0^{(j)}q_0^{(j)}+Blp_0^{(j)}]}E_{Alp_0^{(j)}}T_{-nq_0^{(j)}}\Phi_j(x)\Big\}_{\underset{j\in \{1,2,\dots,N\}}{n,l\in \mathbb{Z}}},
\end{align*}
which is a Gabor system in  $L^2(\mathbb{R},dx)$.

Now we are ready to give necessary and sufficient conditions  for the existence of Gabor frames  for  the Weyl--Heisenberg group.
\begin{thm}\label{w1}
Let $\{\Phi_j\}_{j\in \{1,2,\dots,N\}} \subset L^2(\mathbb{R},dx)$ be a finite collection of non-zero functions such that support of each $\widehat{\Phi_j}$ is contained in an interval of length $\lambda$. For each $j\in \{1,2,\dots,N\}$, let $p_0^{(j)}$ and $q_0^{(j)}$ be real numbers such that
 $|p_0^{(j)}q_0^{(j)}|<1$ and $q_0^{(j)}=\frac{1}{\lambda}$. Then, $\big\{\Phi_{(n,l,j)}^{A,B}\big\}_{n,l\in \mathbb{Z}, j\in \{1,2,\dots,N\}}$ is a Gabor frame for $L^2(\mathbb{R},dx)$ with frame bounds $\alpha_o$,  $\beta_o$ if and only if
	\begin{align*}
\alpha_o \leq \lambda \sum_{j=1}^{N}\sum_{l\in \mathbb{Z}} |\widehat{\Phi_j}(\gamma - Ap_0^{(j)}l)|^2 \leq \beta_o \ \  \text{for almost all} \ \gamma\in \mathbb{R}.
\end{align*}
\end{thm}
For the proof of Theorem \ref{w1}, first we state and prove the  following result which gives necessary and sufficient conditions for nonstationary frames of translates in terms of a series of the Fourier transforms of window functions. Frames of translates (or frames for shift-invariant subspaces) first considered in \cite{ron} by Ron and Shen.

\begin{thm}\label{t1}
For each  $l \in \mathbb{Z}$, let  $\phi_l$ be a  non-zero function in $L^2(\mathbb{R},dx)$ such that support of each $\widehat{\phi_l}$ is contained in an interval of length $\lambda$ and $q^{(l)}=\frac{1}{\lambda}$. Then,  $\{T_{nq^{(l)}}\phi_l\}_{n,l \in \mathbb{Z}}$ is a  nonstationary frame of translates for $L^2(\mathbb{R},dx)$ with frame bounds $\alpha$, $\beta$ if and only if
	\begin{align}\label{t2}
	\frac{\alpha}{\lambda} \leq\sum_{l \in \mathbb{Z}}|\widehat{\phi_l}(\gamma)|^2 \leq \frac{\beta}{\lambda} \ \   \text{for almost all} \ \gamma\in \mathbb{R}.
	\end{align}
\end{thm}

\proof Assume first that \eqref{t2} holds. Then, for any $f \in L^2(\mathbb{R},dx)$, we have
\begin{align}\label{neweq1}
\sum_{n,l\in \mathbb{Z}}|\langle f, T_{nq^{(l)}}\phi_l \rangle|^2 &=\sum_{n,l \in \mathbb{Z}}|\langle \widehat{f}  |  \widehat{T_{nq^{(l)}}\phi_l} \rangle|^2 \notag\\
&=\sum_{n,l \in \mathbb{Z}}|\langle \widehat{f}| E_{-nq^{(l)}}\widehat{\phi_l} \rangle|^2 \notag\\
&=\sum_{n,l \in \mathbb{Z}} \bigg|\int_{\mathbb{R}}\widehat{f}(\gamma)E_{nq^{(l)}}\overline{\widehat{\phi_l}(\gamma)} d\gamma\bigg|^2 \notag\\
&=\sum_{n,l\in \mathbb{Z}} \bigg|\int_{\mathbb{R}}\widehat{f}(\gamma)e^{2\pi inq^{(l)}\gamma}\overline{\widehat{\phi_l}(\gamma)}d\gamma\bigg|^2 \notag\\
&=\sum_{n,l\in \mathbb{Z}} \bigg|\int_{\mathbb{R}}\widehat{f}(\gamma)e^{2\pi in\frac{\gamma}{\lambda}}\overline{\widehat{\phi_l}(\gamma)}d\gamma\bigg|^2.
\end{align}
Using Theorem \ref{thmpp1} in \eqref{neweq1}, we arrive at
\begin{align}\label{eq3.2}
\sum_{n,l\in \mathbb{Z}}|\langle f, T_{nq^{(l)}}\phi_l \rangle|^2 &=\lambda \sum_{l \in \mathbb{Z}}\int_{\mathbb{R}}\bigg|\widehat{f}(\gamma)\overline{\widehat{\phi_l}(\gamma)}\bigg|^2 d\gamma \notag\\
&=\lambda \sum_{l \in \mathbb{Z}}\int_{\mathbb{R}}\big|\widehat{f}(\gamma)\big|^2 \big|\overline{\widehat{\phi_l}(\gamma)}\big|^2 d\gamma \notag\\
&=\lambda \int_{\mathbb{R}}\big|\widehat{f}(\gamma)\big|^2 \sum_{l \in \mathbb{Z}} \big|\overline{\widehat{\phi_l}(\gamma)}\big|^2 d\gamma.
\end{align}
Thus, by Ineq. \eqref{t2}, we have
\begin{align*}
\alpha \int_{\mathbb{R}}\big|\widehat{f}(\gamma)\big|^2 d\gamma \leq \sum_{n,l\in \mathbb{Z}}|\langle f,T_{nq^{(l)}}\phi_l \rangle|^2 \leq \beta \int_{\mathbb{R}}\big|\widehat{f}(\gamma)\big|^2 d\gamma,
\end{align*}
or
\begin{align*}
\alpha \|\hat{f}\|^2 \leq \sum_{n,l\in \mathbb{Z}}|\langle f,T_{nq^{(l)}}\phi_l \rangle|^2 \leq \beta \|\hat{f}\|^2.
\end{align*}
That is
\begin{align*}
\alpha \|f\|^2 \leq \sum_{n,l\in \mathbb{Z}}|\langle f,T_{nq^{(l)}}\phi_l \rangle|^2 \leq \beta \|f\|^2 \ \ \text{for all} \  f \in L^2(\mathbb{R},dx).
\end{align*}
Hence, $\{T_{nq^{(l)}}\phi_l\}_{n,l \in \mathbb{Z}}$ is a  nonstationary frame of translates for $L^2(\mathbb{R},dx)$ with the desired frame bounds.

To prove the opposite implication, let $\{T_{nq^{(l)}}\phi_l\}_{n,l \in \mathbb{Z}}$ be  a  frame of translates for $L^2(\mathbb{R},dx)$ with frame bounds $\alpha, \beta$. Then, using \eqref{eq3.2}, for all $f \in L^2(\mathbb{R}, dx)$, we have
\begin{align*}
\alpha \|f\|^2 \leq \lambda \int_{\mathbb{R}}\big|\widehat{f}(\gamma)\big|^2\sum_{l \in \mathbb{Z}}\big|\overline{\widehat{\phi_l}(\gamma)}\big|^2 d\gamma \leq \beta \|f\|^2,
\end{align*}
which entails
\begin{align*}
\alpha \int_{\mathbb{R}}\big|\widehat{f}(\gamma)\big|^2d\gamma \leq \lambda \int_{\mathbb{R}}\big|\widehat{f}(\gamma)\big|^2\sum_{l \in \mathbb{Z}}\big|\overline{\widehat{\phi_l}(\gamma)}\big|^2 d\gamma \leq \beta \int_{\mathbb{R}}\big|\widehat{f}(\gamma)\big|^2d\gamma
\end{align*}
for all $f \in L^2(\mathbb{R}, dx)$. This completes the proof.
\endproof
\begin{rem}
Lemma 10.1.1 of \cite{ole} is a particular case of  Theorem $\ref{t1}$.
\end{rem}
The following example justifies Theorem $\ref{t1}$.

\begin{exa}\label{et}
Let	$\lambda=\frac{1}{2}$. For each  $l \in \mathbb{Z}$, let $q^{(l)}=2$ and  define functions $\phi_l \in L^2(\mathbb{R},dx)$ as follows.
	\begin{align*}
	\widehat{\phi_l}(x)=\begin{cases}
	2, \quad & \text{if} \ x \in \big(\frac{l}{2},\frac{l+1}{2}\big],\\
	0, \quad & \text{elsewhere}.
	\end{cases}
	\end{align*}
Then, $\sum\limits_{l \in \mathbb{Z}} |\widehat{\phi_l}(\gamma)|^2=2$ for almost all $\gamma \in \mathbb{R}$. Thus, by Theorem ($\ref{t1}$), $\{T_{nq^{(l)}}\phi_l\}_{n,l \in \mathbb{Z}}$ is a nonstationary Parseval  frame of translates for $ L^2(\mathbb{R},dx)$.
\end{exa}

\textbf{Proof of Theorem \ref{w1}:}
First we see that for $n,l\in \mathbb{Z}$ and $j\in \{1,2,\dots,N\}$,
\begin{align}\label{wh1}
\Phi_{(n,l,j)}^{A,B}(x) = e^{2\pi i[(1/2)Anlp_0^{(j)}q_0^{(j)}+Blp_0^{(j)}]}E_{Alp_0^{(j)}}T_{-nq_0^{(j)}}\Phi_j(x),
\end{align}
and
\begin{align}\label{wh2}
E_{Alp_0^{(j)}}T_{-nq_0^{(j)}}\Phi_j=e^{-2\pi iAnlp_0^{(j)}q_0^{(j)}}T_{-nq_0^{(j)}}E_{Alp_0^{(j)}}\Phi_j.
\end{align}
Therefore, for $f\in L^2(\mathbb{R},dx)$, using (\ref{wh1}) and (\ref{wh2}), we have
\begin{align*}
\sum_{j=1}^N\sum_{n,l \in \mathbb{Z}} |\langle f, \Phi_{(n,l,j)}^{A,B}\rangle|^2 &= \sum_{j=1}^N\sum_{n,l \in \mathbb{Z}}  |\langle f,  E_{Alp_0^{(j)}}T_{-nq_0^{(j)}}\Phi_j \rangle |^2\\
&= \sum_{j=1}^N\sum_{n,l \in \mathbb{Z}}  |\langle f,  T_{-nq_0^{(j)}}E_{Alp_0^{(j)}}\Phi_j |^2.
\end{align*}
Thus, $ \{\Phi_{(n,l,j)}^{A,B}\}_{n, l \in \mathbb{Z} \atop j\in \{1,2,\dots,N\}}$ is a frame for $L^2(\mathbb{R},dx)$ if and only if   $\{T_{-nq_0^{(j)}}E_{Alp_0^{(j)}}\Phi_j\}_{n,l \in \mathbb{Z} \atop j\in \{1,2,\dots,N\}}$ is a frame for $L^2(\mathbb{R},dx)$. The result now follows  directly from Theorem $\ref{t1}$.
\endproof

The following example illustrates Theorem \ref{w1}.
\begin{exa}\label{ew}
	Let $A=1$,  $B=0$, $\lambda=2$ and define functions $\Phi_1, \Phi_2 \in L^2(\mathbb{R},dx)$ as follows.
	\begin{align*}
	\widehat{\Phi_1}(\gamma)=\begin{cases}
	1+\gamma, & \text{if}\  \gamma \in (0,1];\\
	\gamma, & \text{if} \ \gamma \in (1,2];\\
	0, &\text{otherwise},
	\end{cases}
	\end{align*}
	and
	\begin{align*}
	\widehat{\Phi_2}(\gamma)=\begin{cases}
	1+\gamma, & \text{if}\  \gamma \in (0,1];\\
	\frac{\gamma}{2}, & \text{if} \ \gamma \in (1,2];\\
	0, &\text{otherwise}.
	\end{cases}
	\end{align*}
	Choose $p_0^{(1)}=p_0^{(2)}=1$. Then  one can easily observe that
	\begin{align*}
	2 \leq \sum_{l \in \mathbb{Z}}|\widehat{\Phi_1}(\gamma-A p_0^{(1)}l)|^2 \leq 8 \ \text{for almost all}\ \gamma \in \mathbb{R};
	\end{align*}
	and
	\begin{align*}
	\frac{5}{4} \leq \sum_{l \in \mathbb{Z}}|\widehat{\Phi_2}(\gamma-A p_0^{(2)}l)|^2 \leq 5 \ \text{for almost all}\ \gamma \in \mathbb{R}.
	\end{align*}
Therefore, we have
\begin{align*}
\frac{13}{2} \leq \lambda\sum_{j=\{1,2\}}\sum_{l \in \mathbb{Z}}|\widehat{\Phi_j}(\gamma-A p_0^{(j)}l)|^2 \leq 26 \ \text{for almost all}\ \gamma \in \mathbb{R}.
\end{align*}
Hence, by Theorem \ref{w1},	$\big\{\Phi_{(n,l,j)}^{1,0}\big\}_{n,l\in \mathbb{Z}, j\in \{1,2\}}$ is a Gabor frame for $L^2(\mathbb{R},dx)$ with frame bounds $\alpha_0=\frac{13}{2}$ and $\beta_0=26$.
\end{exa}

Next,  we discuss an interplay between modulation and translation parameters in Gabor frames for  the Weyl--Heisenberg group. First we observe that if for each $l \in \mathbb{Z}$,  $\phi_l$ is  a  non-zero function in $L^2(\mathbb{R},dx)$ such that support of each $\widehat{\phi_l}$ is contained in an interval of length $\lambda$ and $q^{(l)}=\frac{1}{\lambda}$. Then, using the same approach as in Theorem $\ref{t1}$, we arrive at
\begin{align*}
\sum_{n,l\in \mathbb{Z}}\big|\big\langle f, T_{\frac{n}{q^{(l)}}}\phi_l \big\rangle\big|^2 &=\frac{1}{\lambda} \int_{\mathbb{R}}\big|\widehat{f}(\gamma)\big|^2 \sum_{l \in \mathbb{Z}} \big|\overline{\widehat{\phi_l}(\gamma)}\big|^2 d\gamma\\
&=\frac{1}{\lambda^2}\sum_{n,l\in \mathbb{Z}}|\langle f, T_{nq^{(l)}}\phi_l \rangle|^2.
\end{align*}
Thus, we have the following result.
\begin{prop}\label{whr}
For $l \in \mathbb{Z}$,  $\phi_l$ is  a  non-zero function in $L^2(\mathbb{R},dx)$ such that support of each $\widehat{\phi_l}$ is contained in an interval of length $\lambda$ and $q^{(l)}=\frac{1}{\lambda}$. Then, $\{T_{nq^{(l)}}\phi_l\}_{n,l \in \mathbb{Z}}$ is a  nonstationary frame of translates for $L^2(\mathbb{R},dx)$ with frame bounds $\alpha$, $\beta$ if and only if $\Big\{T_{\frac{n}{q^{(l)}}}\phi_l\Big\}_{n,l \in \mathbb{Z}}$ is a  nonstationary frame of translates for $L^2(\mathbb{R},dx)$ with frame bounds $\frac{\alpha}{\lambda^2}$, $\frac{\beta}{\lambda^2}$.
\end{prop}

Now, we are ready to give an  interplay between modulation and translation parameters  for  Gabor frames associated with the Weyl--Heisenberg group.
\begin{thm}\label{3.77}
For $j \in \{1,2,\dots,\mathbb{N}\}$, let  $\Phi_j$ be  non-zero functions in $L^2(\mathbb{R},dx)$ such that support of each $\widehat{\Phi_j}$ is contained in an interval of length $\lambda$ and $q_0^{(j)}=\frac{1}{\lambda}$. Then,  $\big\{\Phi_{(n,l,j)}^{A,B}\big\}_{n,l \in \mathbb{Z},j \in \{1,2,\dots,N\}}$ is a Gabor frame  for $L^2(\mathbb{R},dx)$ with frame bounds $\alpha, \beta$ if and only if	
$\Big\{ E_{Alp_0^{(j)}}T_{\frac{n}{q_0^{(j)}}}\Phi_j\Big\}_{n,l \in \mathbb{Z},j \in \{1,2,\dots,N\}}$ is a Gabor frame for $L^2(\mathbb{R},dx)$ with frame bounds $\frac{\alpha}{\lambda^2}$, $\frac{\beta}{\lambda^2}$.
\end{thm}

\proof
Since
\begin{align}\label{wif1}
E_{Alp_0^{(j)}}T_{-nq_0^{(j)}}\Phi_j=e^{-2\pi iAnlp_0^{(j)}q_0^{(j)}}T_{-nq_0^{(j)}}E_{Alp_0^{(j)}}\Phi_j  \ \text{for all} \ n, l \in \mathbb{Z} \ \text{and} \ j \in \{1,2,\dots,N\}.
\end{align}
Therefore,  from Eq. (\ref{wif1}) and Proposition \ref{whr}, we can say that
\begin{align*}
&\big\{\Phi_{(n,l,j)}^{A,B}\big\}_{n,l \in \mathbb{Z},j \in \{1,2,\dots,N\}} \ \text{is a Gabor frame  for} \ L^2(\mathbb{R},dx) \  \text{with frame bounds} \ \alpha, \beta. \\
& \Leftrightarrow \big\{T_{-nq_0^{(j)}}\Phi_{(l,j)}^{A,B}\big\}_{n,l \in \mathbb{Z},j \in \{1,2,\dots,N\}} \  \text{is a nonstationary frame of translates for} \ L^2(\mathbb{R},dx) \ \text{with frame}\\
& \quad \ \text{bounds} \  \alpha, \beta, \ \text{where} \ \Phi_{l,j}=E_{Alp_0^{(j)}}\Phi_j.\\
& \Leftrightarrow  \Big\{T_{\frac{-n}{q_0^{(j)}}}\Phi_{(l,j)}^{A,B}\Big\}_{n,l \in \mathbb{Z},j \in \{1,2,\dots,N\}} \text{is a nonstationary frame of translates for}  \ L^2(\mathbb{R},dx) \ \text{with frame} \\
& \quad  \ \text{bounds} \ \frac{\alpha}{\lambda^2}, \frac{\beta}{\lambda^2}.\\
&\Leftrightarrow  \Big\{E_{Alp_0^{(j)}}T_{\frac{n}{q_0^{(j)}}}\Phi_j\Big\}_{n,l \in \mathbb{Z},j \in \{1,2,\dots,N\}} \ \text{is a Gabor frame for} \ L^2(\mathbb{R},dx)  \ \text{with frame bounds} \ \frac{\alpha}{\lambda^2}, \frac{\beta}{\lambda^2}.
\end{align*}
This concludes the proof.
\endproof
\begin{lem}\cite[Lemma 5.3.3]{ole}\label{l1}
If $\{\psi_k\}_{k \in \mathbb{Z}}$ is a frame sequence with bounds $\alpha$ and  $\beta$ and $U \colon \mathcal{H} \rightarrow \mathcal{H}$ is a unitary operator, then $\{U\psi_k\}_{k \in \mathbb{Z}}$ is a frame sequence with  bounds $\alpha$ and  $\beta$.
\end{lem}

\begin{prop}
	 Assume that $\{\Phi_{(n,l,j)}^{A,B}\}_{n,l \in \mathbb{Z},j \in \{1,2,\dots,N\}}$ is a frame sequence  with frame bounds $\alpha$, $\beta$. Given $a>0$, let $\phi_{j_a}:=D_a\phi_j$. Then, $\{\Phi_{(na,\frac{l}{a},j_a)}^{A,B}\}_{n,l \in \mathbb{Z},j \in \{1,2,\dots,N\}}$ is a  frame sequence  with  bounds $\alpha$ and  $\beta$.
\end{prop}

\proof
Using the commutator relations of Lemma \ref{p1}, we have
\begin{align*}
\Phi_{(na,\frac{l}{a},j_a)}^{A,B}(x)&= e^{2\pi i[(1/2)Anlp_0^{(j)}q_0^{(j)}+Blp_0^{(j)}]}E_{A\frac{l}{a}p_0^{(j)}}T_{-naq_0^{(j)}}D_a\Phi_j(x)\\
&= e^{2\pi i[(1/2)Anlp_0^{(j)}q_0^{(j)}+Blp_0^{(j)}]}E_{A\frac{l}{a}p_0^{(j)}}D_a T_{-nq_0^{(j)}}\Phi_j(x)\\
&= e^{2\pi i[(1/2)Anlp_0^{(j)}q_0^{(j)}+Blp_0^{(j)}]}D_a E_{Alp_0^{(j)}}T_{-nq_0^{(j)}}\Phi_j(x)\\
&= D_a \Phi_{(n,l,j)}^{A,B}(x).
\end{align*}
The proof now follows from Lemma $\ref{l1}$.
\endproof
It is proved in \cite[Lemma 9.3.2]{ole} that a  function in $L^2(\mathbb{R})$ belongs to the  closure of the span of a frame sequence of translates if and only if its Fourier transform can be expressed in terms of the Fourier transform of the window functions. The next result generalises \cite[Lemma 9.3.2]{ole} for frame sequences  with the Gabor structure associated with the Weyl--Heisenberg group.
\begin{thm}\label{thm3.10x}
	Assume that $\{\Phi_{(n,l,j)}^{A,B}\}_{n,l \in \mathbb{Z},j \in \{1,2,\dots,N\}}$ is a  frame sequence  in $L^2(\mathbb{R},dx)$. Then, a function $f \in L^2(\mathbb{R},dx)$ belongs to  $\overline{\text{span}}\{\Phi_{(n,l,j)}^{A,B}\}_{n,l \in \mathbb{Z},j \in \{1,2,\dots,N\}}$ if and only if there exists a sequence $\{F_{l,j}\}_{l \in \mathbb{Z},j\in\{1,2,\dots,N\}}$ of $1$-periodic functions such that
	\begin{align*}
	\widehat{f}(\gamma)=\sum_{j=1}^{N}\sum_{l \in \mathbb{Z}}F_{l,j}\widehat{\Phi_j}(\gamma-Alp_0^{(j)}),
	\end{align*}
	where restriction of each $F_{l,j}$ to $[0,1)$ belongs to $L^2([0,1),dx)$.
\end{thm}
\proof
A function $f\in L^2(\mathbb{R},dx)$ belongs to $\overline{\text{span}}\{\Phi_{(n,l,j)}^{A,B}\}_{n,l \in \mathbb{Z},j \in \{1,2,\dots,N\}}$ if and only if there exists a sequence $\{c_{n,l,j}\}_{n,l\in \mathbb{Z},j\in \{1,2,\dots,N\}}\in \ell^2(\mathbb{Z} \times \mathbb{Z}\times \{1,2,\dots,N\})$ such that
\begin{align*}
f=\sum_{j=1}^{N}\sum_{n,l\in \mathbb{Z}} c_{n,l,j}\Phi_{(n,l,j)}^{A,B}.
\end{align*}
Therefore
\begin{align*}
\widehat{f}(\gamma)&=\sum_{j=1}^{N}\sum_{n,l\in \mathbb{Z}} c_{n,l,j}e^{2\pi i[(1/2)Anlp_0^{(j)}q_0^{(j)}+Blp_0^{(j)}]}e^{-2\pi iAnlp_0^{(j)}q_0^{(j)}}E_{nq_0^{(j)}}T_{Alp_0^{(j)}}\widehat{\Phi_j}(\gamma)\\
&= \sum_{j=1}^{N}\sum_{n,l\in \mathbb{Z}} c_{n,l,j}e^{2\pi i[(-1/2)Anlp_0^{(j)}q_0^{(j)}+Blp_0^{(j)}]}e^{2\pi i nq_0^{(j)}\gamma}\widehat{\Phi_j}(\gamma-Alp_0^{(j)})\\
&=\sum_{l \in \mathbb{Z}} F_l \widehat{\Phi_l}(\gamma),
\end{align*}
where $F_l=\sum\limits_{j=1}^{N}\sum\limits_{n \in \mathbb{Z}} c_{n,l,j}e^{2\pi i[(-1/2)Anlp_0^{(j)}q_0^{(j)}+Blp_0^{(j)}]}e^{2\pi i nq_0^{(j)}\gamma}$. This concludes the result.
\endproof

\section{Frames for the Extended Affine Group}\label{Sec4}
The affine group $A$ realized as the upper half plane
\begin{align*}
\mathcal{A}=\{(\alpha,\beta) \colon \alpha \in \mathbb{R}^+, \beta \in \mathbb{R}\}.
\end{align*}
The multiplication law in $A$ is given by
\begin{align*}
(\alpha,\beta)(x, y)=(\alpha x, \alpha y +\beta)
\end{align*}
for all points $(\alpha, \beta)$, $(x,y)\in A$.

Let $V$ denote the subspace of $L^2(\mathbb{R},dx)$ consisting of all the functions whose Fourier transform vanish on the negative half of $\mathbb{R}$.
Recall that for any $d \in \mathbb{R}^*$, the map
\begin{align*}
\rho^d \colon A &\rightarrow \mathcal{U}(V)\\
(\rho^d(\alpha,\beta)f)(x) &=\frac{1}{\sqrt{\alpha}}f\Big(\frac{x+d\beta}{\alpha}\Big), \ x \in \mathbb{R}.
\end{align*}

The extended affine group, denoted by $EA$, is the direct product of the affine group $A$ with the additive group of real numbers. That is, $EA=A \oplus \mathbb{R}$. The product of $(\alpha,\beta,\gamma)$ and $(x,y,z)$ in $EA$ is given by
\begin{align*}
(\alpha,\beta,\gamma)(x,y,z)=(\alpha x,\alpha y +\beta, \gamma +z).
\end{align*}

For $c \in \mathbb{R}$, let $\chi_c(x)=e^{icx}$ be a unitary character of $\mathbb{R}$. For every $d\in \mathbb{R}^*, c \in \mathbb{R}$, the unitary irreducible representation of $EA$ acting on $V$ is the map
\begin{align*}
\beth^{c,d}=\rho^d \otimes \chi_c \colon EA &\rightarrow \mathcal{U}(V)\\
(\beth^{c,d}(\alpha,\beta,\gamma)f)(x)&=e^{ic\gamma}\frac{1}{\alpha}f\Big(\frac{x+d\beta}{\alpha}\Big).
\end{align*}

We now study wavelet frames associated with the extended affine group. We fix a representation $\Theta^{A,B}$ of the Weyl--Heisenberg group on $L^2(\mathbb{R})$. For $j \in \{1,2,\dots,N\}$, where $N$ is some strictly positive number, consider the discrete subset of $EA$ that is defined by $EA_{q_0^{(j)},p_0^{(j)}}=\{(\alpha_l,\beta_{nl},\gamma_{nl}) \colon n,l\in \mathbb{Z}\},$ where
\begin{align*}
\alpha_l&=e^{-lq_0^{(j)}},\\ \beta_{nl}&=\frac{-nl(q_0^{(j)})^2}{d\hspace{.5mm} \mathrm{ln}\hspace{.5mm}(\alpha_l)}=\frac{nq_0^{(j)}}{d}, \\ \gamma_{nl}&=\beta_{nl}\frac{\mathrm{ln}(\alpha_l)}{\alpha_l-1}.
\end{align*}

For $\Phi_j \in L^2(\mathbb{R},dx)\backslash\{0\}$ and  $j\in \{1,2,\dots,N\}$, define
\begin{align*}
\Phi_{(n,l,j)}^{c,d}=\beth^{c,d}(\alpha_l, \beta_{nl}, \gamma_{nl})\Phi_j,
\end{align*}
where $n,l \in \mathbb{Z}$. Thus,
\begin{align}\label{ea1}
\Phi_{(n,l,j)}^{c,d}(x)&=e^{-\frac{i c \mathrm{ln}(q_0^{(j)})^2}{d(e^{-lq_0^{(j)}}-1)}} \frac{1}{\sqrt{e^{-lq_0^{(j)}}}} \Phi_j(\frac{x+nq_0^{(j)}}{e^{-l q_0^{(j)}}})\notag\\
&=e^{-\frac{i c \mathrm{ln}(q_0^{(j)})^2}{d(e^{-lq_0^{(j)}}-1)}}D_{e^{-lq_0^{(j)}}}T_{\frac{-nq_0^{(j)}}{e^{-lq_0^{(j)}}}}\Phi_j(x), \ \text{for} \ n,l \in \mathbb{Z}, j \in \{1,2,\dots,N\},
\end{align}
which is equivalent to wavelet system. For more technical details about frames associated with extended affine group, we refer to \cite{SB}.

The following result gives the existence of wavelet frame associated with the  extended affine group.

\begin{thm}\label{th4.8e}
	Under the assumptions of Theorem $\ref{w1}$, the sequence $\{\Phi_{(n,l,j)}^{c,d}\}_{\underset{j \in \{1,2,\dots,N\}}{n,l\in \mathbb{Z}}}$ is a wavelet frame for $L^2(\mathbb{R},dx)$ with bounds $\alpha_o$, $\beta_o$ if and only if
	\begin{align*}
\alpha_o \leq \lambda \sum_{j=1}^{N}\sum_{l \in \mathbb{Z}}\frac{1}{e^{lq_0^{(j)}}}\bigg|\widehat{\Phi_j}\bigg(\frac{x}{e^{lq_0^{(j)}}}\bigg)\bigg|^2 \leq \beta_o \ \text{for almost all} \ \gamma \in \mathbb{R}.
	\end{align*}
\end{thm}
\proof
For all $n,l\in \mathbb{Z}$ and $j\in \{1,2,\dots,N\}$, we have
\begin{align}\label{ea2}
D_{e^{-lq_0^{(j)}}}T_{\frac{-nq_0^{(j)}}{e^{-lq_0^{(j)}}}}\Phi_j= T_{-nq_0^{(j)}}D_{e^{-lq_0^{(j)}}}\Phi_j.
\end{align}
Using (\ref{ea1}) and (\ref{ea2}), for $f\in L^2(\mathbb{R},dx)$, we have
\begin{align*}
\sum_{j=1}^N\sum_{n,l \in \mathbb{Z}} |\langle f, \Phi_{(n,l,j)}^{c,d}\rangle|^2 &= \sum_{j=1}^N\sum_{n,l \in \mathbb{Z}}  |\langle f, D_{e^{-lq_0^{(j)}}}T_{\frac{-nq_0^{(j)}}{e^{-lq_0^{(j)}}}}\Phi_j(x) \rangle |^2\\
&= \sum_{j=1}^N\sum_{n,l \in \mathbb{Z}}  |\langle f, T_{-nq_0^{(j)}}D_{e^{-lq_0^{(j)}}}\Phi_j |^2.
\end{align*}
Thus, $ \big\{\Phi_{(n,l,j)}^{c,d}\}_{n,l \in \mathbb{Z} \atop  j\in \{1,2,\dots,N\}}$ is a frame for  $L^2(\mathbb{R},dx)$ if and only if the   $\big\{T_{-nq_0^{(j)}}D_{e^{-lq_0^{(j)}}}\Phi_j\big\}_{n,l \in \mathbb{Z} \atop j\in \{1,2,\dots,N\}}$ is a frame for $L^2(\mathbb{R},dx)$. The result now directly follows from Theorem $\ref{t1}$.
\endproof
\section{The Canonical Dual of Nonstationary Frames of Translates}\label{Sec5}
Recently, the structure of the canonical dual of different types of frames studied by many authors, see \cite{ CS, DVI, CK} and references therein.  Chui and Shi showed in \cite{CS} that the canonical dual of a wavelet frame for $L^2(\mathbb{R})$ need not have a wavelet structure. The authors of \cite{DVI} proved that  the canonical dual of a  discrete wavelet frame for the unitary space $\mathbb{C}^N$ has the same structure. Further, the canonical dual of a discrete Gabor frame has the same structure, see Proposition 6.1 of  \cite{CK}. In this section, we show that the canonical dual of  nonstationary frames of translates has the same structure.
\begin{thm}\label{th5.1}
Let $\{T_{nq^{(l)}}\phi_l\}_{n,l \in \mathbb{Z}}$ be a nonstationary frame of translates for $L^2(\mathbb{R},dx)$ with frame operator $\mathcal{S}$. Then, the canonical dual frame of $\{T_{nq^{(l)}}\phi_l\}_{n,l \in \mathbb{Z}}$ is $\{T_{nq^{(l)}}\mathcal{S}^{-1}\phi_l\}_{n,l \in \mathbb{Z}}$.
\end{thm}

\proof
First we show that frame operator $\mathcal{S}$ commutes with translation operator.
For any $n',l'\in \mathbb{Z}$ and $\psi_{l'}\in L^2(\mathbb{R},dx)$, we compute
\begin{align*}
T_{n'q^{(l')}}\mathcal{S}\psi_{l'}&=T_{n'q^{(l')}} \sum_{n,l \in \mathbb{Z}}\langle \psi_{l'} , T_{nq^{(l)}}\phi_l \rangle T_{nq^{(l)}}\phi_l \\
&=\sum_{n,l \in \mathbb{Z}} \langle \psi_{l'} ,T_{nq^{(l)}}\phi_l \rangle T_{n'q^{(l')}} T_{nq^{(l)}}\phi_l\\
&=\sum_{n,l \in \mathbb{Z}} \langle \psi_{l'} ,T_{nq^{(l)}}\phi_l \rangle T_{n'q^{(l')}+nq^{(l)}}\phi_l\\
&=\sum_{n,l \in \mathbb{Z}} \langle \psi_{l'},T_{nq^{(l)}-n'q^{(l')}}\phi_l \rangle T_{nq^{(l)}} \phi_l\\
&= \sum_{n,l \in \mathbb{Z}} \langle \psi_{l'}, T_{-n'q^{(l')}}T_{nq^{(l)}}\phi_l \rangle T_{nq^{(l)}} \phi_l \\
&= \sum_{n,l \in \mathbb{Z}} \langle T_{-n'q^{(l')}}\psi_{l'}, T_{nq^{(l)}}\phi_l \rangle T_{nq^{(l)}} \phi_l\\
&=\mathcal{S} T_{-n'q^{(l')}}\psi_{l'}.
\end{align*}
Therefore, the frame operator $\mathcal{S}$ commutes with translation operator.

Now
\begin{align*}
\mathcal{S}^{-1} T_{nq^{(l)}} \phi_l&=(T_{nq^{(l)}}^{-1}\mathcal{S})^{-1} \phi_l\\
&= (T_{-nq^{(l)}}\mathcal{S})^{-1} \phi_l\\
&=  (\mathcal{S}T_{-nq^{(l)}})^{-1} \phi_l \\
&= T_{-nq^{(l)}} \mathcal{S}^{-1}\phi_l, \ n, l \in \mathbb{Z}.
\end{align*}
Hence, the canonical dual frame of $\{T_{nq^{(l)}}\phi_l\}_{n,l \in \mathbb{Z}}$ is $\big\{T_{nq^{(l)}}\mathcal{S}^{-1}\phi_l\big\}_{n,l \in \mathbb{Z}}$. This completes the proof.
\endproof
\section{Nonstationary Frames and Riesz Bases of Translates}\label{Sec6}
In   \cite{Casa},  Casazza and  Christensen introduced a new method to approximate the inverse of the frame operator using finite subsets of the frame. In their study, they  also consider Gabor frames and frames consisting of translates of a single function. In this section, we consider the  technique given in \cite{Casa} for nonstationary frames of translates.  Consider a nonstationary frames of translates $\{T_{nq^{(l)}}\phi_l\}_{n,l \in \mathbb{Z}}$ of $L^2(\mathbb{R})$ with frame operator $\mathcal{S}$, where $q^{(l)} \in \mathbb{R}^{+}$ for each $l \in \mathbb{Z}$. We recall that every finite collection of vectors  in a Hilbert space $\mathcal{H}$ is a frame for its span, see \cite{ole}.  Then, the family $\{T_{nq^{(l)}}\phi_l\}_{n\in \mathbb{F}_1, l\in \mathbb{F}_2}$ is a frame for   $\mathcal{V}=\text{span}\{T_{nq^{(l)}}\phi_l\}_{n\in \mathbb{F}_1, l\in \mathbb{F}_2}$, where $\mathbb{F}_1=\{-s,\dots,s\}, \mathbb{F}_2=\{-t,\dots,t\}$ for $s,t \in \mathbb{N}$. Denote its frame operator by $\mathcal{S}_{s,t}$. Then,  $\mathcal{S}_{s,t} \colon \mathcal{V} \rightarrow \mathcal{V}$ is given by
\begin{align*}
\mathcal{S}_{s,t}f =\sum_{n \in \mathbb{F}_1}\sum_{l \in \mathbb{F}_2} \langle f , T_{nq^{(l)}}\phi_l\rangle T_{nq^{(l)}}\phi_l,
\end{align*}
and its frame decomposition is given by  $f=\sum_{n \in \mathbb{F}_1}\sum_{l \in \mathbb{F}_2} \langle f , \mathcal{S}_{s,t}^{-1}T_{nq^{(l)}}\phi_l\rangle T_{nq^{(l)}}\phi_l$, $f\in \mathcal{V}$. Since a nonstationary frames of translates are multivariate frames, it would be interesting to approximate the inverse of frame operator (of a multivariate frame) in some sense.
In this direction, the following result gives a necessary and sufficient condition for approximation of inverse of the frame operator of nonstationary  frames of translates in the weak sense.

\begin{thm}\label{r1}
For each $l \in \mathbb{Z}$, let $q^{(l)} \in \mathbb{R}^{+}$. Suppose $\{T_{nq^{(l)}}\phi_l\}_{n,l \in \mathbb{Z}}$ is  a nonstationary frame of translates  for $L^2(\mathbb{R})$ with upper bound $\beta$. Then, for all $f \in L^2(\mathbb{R})$  and for all  $n, l \in \mathbb{Z}$,
	\begin{align}\label{s1}
	\big\langle f , \mathcal{S}_{s,t}^{-1}T_{nq^{(l)}}\phi_l \big\rangle \rightarrow \big \langle f, \mathcal{S}^{-1}T_{nq^{(l)}}\phi_l \big\rangle \ \text{as} \ s, t  \rightarrow \infty,
	\end{align}
 if and only if, for all  $i$, $j \in \mathbb{N}$, there exists  $c_{i,j}\in \mathbb{R}$ such that
\begin{align}\label{s2}
 \|T_{iq^{(j)}}\mathcal{S}_{s,t}^{-1}\phi_j\|= \|\mathcal{S}_{s,t}^{-1}\phi_j\| \leq c_{i,j} \ \text{for all}\  \ s\geq i, \ t \geq j.
\end{align}
\end{thm}
\proof
Suppose first that $(\ref{s2})$ holds. For $i$, $j\in \mathbb{N}$ define functions $\psi_{s,t}$ as follows.
\begin{align*}
\psi_{s,t}= \mathcal{S}_{s,t}^{-1}T_{iq^{(j)}}\phi_j - \mathcal{S}^{-1}T_{iq^{(j)}}\phi_j.
\end{align*}
We need to prove that for all $f \in L^2(\mathbb{R})$,  $\langle f , \psi_{s,t}\rangle \to 0$ as $s,t \to \infty$.

Note that for all $f \in L^2(\mathbb{R})$, we have
\begin{align}\label{s3}
\mathcal{S} f &= \sum_{n,l \in \mathbb{Z}} \langle f , T_{nq^{(l)}}\phi_l\rangle T_{nq^{(l)}}\phi_l \nonumber \\
&= \mathcal{S}_{s,t}f + \sum_{n\in \mathbb{Z}}\sum_{l\in \mathbb{Z} \setminus \mathbb{F}_2}\langle f , T_{nq^{(l)}}\phi_l\rangle T_{nq^{(l)}}\phi_l + \sum_{n \in \mathbb{Z} \setminus \mathbb{F}_1}\sum_{l\in \mathbb{Z}}\langle f , T_{nq^{(l)}}\phi_l\rangle T_{nq^{(l)}}\phi_l.
\end{align}
We will use this to obtain an alternative formula for $\psi_{s,t}$.

Since
\begin{align*}
\mathcal{S}\psi_{s,t}= \mathcal{S}\mathcal{S}_{s,t}^{-1}T_{iq^{(j)}}\phi_j -T_{iq^{(j)}}\phi_j.
\end{align*}
Therefore, using Eq.  $(\ref{s3})$ on $\mathcal{S}_{s,t}^{-1}T_{iq^{(j)}}\phi_j$, we obtain
\begin{align*}
\mathcal{S}\psi_{s,t}&= \mathcal{S}_{s,t}\mathcal{S}_{s,t}^{-1}T_{iq^{(j)}}\phi_j+ \sum_{n \in \mathbb{Z}}\sum_{l \in \mathbb{Z}\setminus\mathbb{F}_2}^{\infty}\langle \mathcal{S}_{s,t}^{-1}T_{iq^{(j)}}\phi_j , T_{nq^{(l)}}\phi_l\rangle T_{nq^{(l)}}\phi_l \\
&  + \sum_{n \in \mathbb{Z}\setminus\mathbb{F}_1}\sum_{l \in \mathbb{Z}}\langle \mathcal{S}_{s,t}^{-1}T_{iq^{(j)}}\phi_j , T_{nq^{(l)}}\phi_l\rangle T_{nq^{(l)}}\phi_l -
 \ T_{iq^{(j)}}\phi_j \\
&= \sum_{n\in \mathbb{Z}}\sum_{l \in \mathbb{Z}\setminus\mathbb{F}_2}^{\infty}\langle \mathcal{S}_{s,t}^{-1}T_{iq^{(j)}}\phi_j , T_{nq^{(l)}}\phi_l\rangle T_{nq^{(l)}}\phi_l + \sum_{n \in \mathbb{Z}\setminus\mathbb{F}_1}\sum_{l \in \mathbb{Z}}\langle \mathcal{S}_{s,t}^{-1}T_{iq^{(j)}}\phi_j , T_{nq^{(l)}}\phi_l\rangle T_{nq^{(l)}}\phi_l.
\end{align*}
It follows that for $s\geq i$ and $t\geq j$,
\begin{align*}
\psi_{s,t} = \sum_{n \in \mathbb{Z}}\sum_{l \in \mathbb{Z}\setminus\mathbb{F}_2}\langle \mathcal{S}_{s,t}^{-1}T_{iq^{(j)}}\phi_j , T_{nq^{(l)}}\phi_l\rangle \mathcal{S}^{-1}T_{nq^{(l)}}\phi_l + \sum_{n \in \mathbb{Z}\setminus\mathbb{F}_1}\sum_{l\in \mathbb{Z}}\langle \mathcal{S}_{s,t}^{-1}T_{iq^{(j)}}\phi_j , T_{nq^{(l)}}\phi_l\rangle \mathcal{S}^{-1}T_{nq^{(l)}}\phi_l.
\end{align*}
Therefore, for all $f \in L^2(\mathbb{R},dx)$, we have
\begin{align}\label{ap1}
|\langle f , \psi_{s,t} \rangle|^2 &=\Big|\sum_{n\in \mathbb{Z}}\sum_{l\in \mathbb{Z}\setminus\mathbb{F}_2}\langle \mathcal{S}_{s,t}^{-1}T_{iq^{(j)}}\phi_j , T_{nq^{(l)}}\phi_l\rangle\langle f ,\mathcal{S}^{-1}T_{nq^{(l)}}\phi_l\rangle \notag \\
&  + \sum_{n\in \mathbb{Z}\setminus\mathbb{F}_1}\sum_{l\in \mathbb{Z}}\langle \mathcal{S}_{s,t}^{-1}T_{iq^{(j)}}\phi_j , T_{nq^{(l)}}\phi_l\rangle \langle f,\mathcal{S}^{-1}T_{nq^{(l)}}\phi_l\rangle \Big|^2 \notag\\
&\leq 2\Big(\Big|\sum_{n\in \mathbb{Z}}\sum_{l\in \mathbb{Z}\setminus\mathbb{F}_2}\langle \mathcal{S}_{s,t}^{-1}T_{iq^{(j)}}\phi_j , T_{nq^{(l)}}\phi_l\rangle\langle f ,\mathcal{S}^{-1}T_{nq^{(l)}}\phi_l\rangle \Big|^2 \notag\\
& + \Big|\sum_{n\in \mathbb{Z}\setminus\mathbb{F}_1}\sum_{l\in \mathbb{Z}}\langle \mathcal{S}_{s,t}^{-1}T_{iq^{(j)}}\phi_j , T_{nq^{(l)}}\phi_l\rangle \langle f,\mathcal{S}^{-1}T_{nq^{(l)}}\phi_l\rangle \Big|^2\Big).
\end{align}
Using Cauchy-Schwartz' inequality in \eqref{ap1}, we compute
\begin{align}\label{ap2}
|\langle f , \psi_{s,t} \rangle|^2 &\leq 2\Big(\sum_{n\in \mathbb{Z}}\sum_{l\in \mathbb{Z}\setminus\mathbb{F}_2}|\langle \mathcal{S}_{s,t}^{-1}T_{iq^{(j)}}\phi_j , T_{nq^{(l)}}\phi_l\rangle|^2 \sum_{n\in \mathbb{Z}}\sum_{l\in \mathbb{Z}\setminus\mathbb{F}_2}|\langle f ,\mathcal{S}^{-1}T_{nq^{(l)}}\phi_l\rangle|^2 + \notag\\
&  +  \sum_{n\in \mathbb{Z}\setminus\mathbb{F}_1}\sum_{l\in \mathbb{Z}} |\langle \mathcal{S}_{s,t}^{-1}T_{iq^{(j)}}\phi_j , T_{nq^{(l)}}\phi_l\rangle|^2\sum_{n\in \mathbb{Z}\setminus\mathbb{F}_1}\sum_{l\in \mathbb{Z}}|\langle f,\mathcal{S}^{-1}T_{nq^{(l)}}\phi_l\rangle |^2\Big) \notag\\
& \leq 2\beta \|\mathcal{S}_{s,t}^{-1}T_{iq^{(j)}}\phi_j \|^2 \bigg(\sum_{n\in \mathbb{Z}}\sum_{l\in \mathbb{Z}\setminus\mathbb{F}_2} |\langle f ,\mathcal{S}^{-1}T_{nq^{(l)}}\phi_l\rangle|^2 + \sum_{n\in \mathbb{Z}\setminus\mathbb{F}_1}\sum_{l\in \mathbb{Z}} |\langle f ,\mathcal{S}^{-1}T_{nq^{(l)}}\phi_l\rangle|^2 \bigg) \notag\\
&= 2\beta \|T_{iq^{(j)}}\mathcal{S}_{s,t}^{-1}\phi_j \|^2 \bigg(\sum_{n\in \mathbb{Z}}\sum_{l\in \mathbb{Z}\setminus\mathbb{F}_2} |\langle f ,\mathcal{S}^{-1}T_{nq^{(l)}}\phi_l\rangle|^2 + \sum_{n\in \mathbb{Z}\setminus\mathbb{F}_1}\sum_{l\in \mathbb{Z}} |\langle f ,\mathcal{S}^{-1}T_{nq^{(l)}}\phi_l\rangle|^2 \bigg) \notag\\
&= 2\beta \|\mathcal{S}_{s,t}^{-1}\phi_j \|^2 \bigg(\sum_{n\in \mathbb{Z}}\sum_{l\in \mathbb{Z}\setminus\mathbb{F}_2} |\langle f ,\mathcal{S}^{-1}T_{nq^{(l)}}\phi_l\rangle|^2 + \sum_{n\in \mathbb{Z}\setminus\mathbb{F}_1}\sum_{l\in \mathbb{Z}} |\langle f ,\mathcal{S}^{-1}T_{nq^{(l)}}\phi_l\rangle|^2 \bigg) \notag\\
&\leq 2 \beta c_{i,j}^2 \bigg(\sum_{n\in \mathbb{Z}}\sum_{l\in \mathbb{Z}\setminus\mathbb{F}_2} |\langle \mathcal{S}^{-1}f ,T_{nq^{(l)}}\phi_l\rangle|^2 + \sum_{n\in \mathbb{Z}\setminus\mathbb{F}_1}\sum_{l\in \mathbb{Z}} |\langle \mathcal{S}^{-1} f ,T_{nq^{(l)}}\phi_l\rangle|^2 \bigg).
\end{align}
Since $\{T_{nq^{(l)}}\phi_l\}_{n,l \in \mathbb{Z}}$ is a frame, we have
\begin{align}\label{ap3}
\sum_{n\in \mathbb{Z}}\sum_{l\in \mathbb{Z}\setminus\mathbb{F}_2} |\langle \mathcal{S}^{-1} f ,T_{nq^{(l)}}\phi_l\rangle|^2 \to 0 \quad  \text{and}  \quad  \sum_{n\in \mathbb{Z}\setminus\mathbb{F}_1}\sum_{l\in \mathbb{Z}} |\langle \mathcal{S}^{-1}f ,T_{nq^{(l)}}\phi_l\rangle|^2 \to 0 \ \text{as} \ s, t\to \infty.
\end{align}
From \eqref{ap2} and \eqref{ap3}, we conclude that
\begin{align*}
\langle f , \psi_{s,t}\rangle \to 0 \ \text{as} \ s,t \to \infty
\end{align*}
for all $f \in L^2(\mathbb{R})$.

To prove the converse,  assume that $(\ref{s1})$ holds. For any $i$, $j \in \mathbb{N}$, consider the functionals
\begin{align*}
 F_{s,t} &\colon L^2(\mathbb{R},dx) \rightarrow \mathbb{C},\\
  F_{s,t} f &= \langle f ,\mathcal{S}_{s,t}^{-1}T_{iq^{(j)}}\phi_j\rangle, \quad s\geq i, \  t\geq j.
\end{align*}
Each $F_{s,t}$ is bounded and linear. By condition $(\ref{s1})$,  the family of functionals $\{F_{s,t}\}_{s\geq i, t\geq j}$ is pointwise convergent. By invoking the Principal of Uniform Boundedness, the sequence $\{\|F_{s,t}\|\}_{s\geq i, t\geq j}$ is bounded. Therefore, there is constant $c_{i,j}>0$ such that
\begin{align*}
c_{i,j}\geq \|F_{s,t}\| &= \|\mathcal{S}_{s,t}^{-1}T_{iq^{(j)}}\phi_j\|
= \|T_{iq^{(j)}}\mathcal{S}_{s,t}^{-1}\phi_j\|
= \|\mathcal{S}_{s,t}^{-1}\phi_j\| \ \text{for all} \ s\geq i, \  t\geq j.
\end{align*}
This completes the proof.
\endproof

\begin{exa}
	Consider the nonstationary  Parseval frame $\{T_{nq^{(l)}}\phi_l\}_{n,l \in \mathbb{Z}}$ given in Example \ref{et}. Then, $\alpha_{Opt}=\beta_{Opt}=1$ are optimal frame bounds for $\{T_{nq^{(l)}}\phi_l\}_{n,l \in \mathbb{Z}}$. Observe that  $\|\phi_l\| = 1$ for every $l \in \mathbb{Z}$. Hence, $\|T_{nq^{(l)}}\phi_l\| =1$; $n,l \in \mathbb{Z}$. For all $s \geq i$, $t \geq j$,
	\begin{align*}
		\|S_{s,t}^{-1}\| = \sup_{\|T_{iq^{(j)}}\phi_j\|=1} \|S_{s,t}^{-1}T_{iq^{(j)}}\phi_j\| \geq \|S_{s,t}^{-1}T_{iq^{(j)}}\phi_j\| = \|T_{iq^{(j)}}S_{s,t}^{-1}\phi_j\|
	\end{align*}
	By Proposition \ref{o1}, $\|S_{s,t}\| = \|S_{s,t}^{-1}\| = 1$ for all $s \geq i$, $t \geq j$. Thus,
	\begin{align*}
		\|T_{iq^{(j)}}S_{s,t}^{-1}\phi_j\| \leq 1 \ \text{for all}\ s \geq i, t \geq j.
	\end{align*}
	Hence, by Theorem \ref{r1}, \eqref{s1} holds.
\end{exa}

For a nonstationary  frame of translates $\{T_{nq^{(l)}}\phi_l\}_{n,l \in \mathbb{Z}}$,  the following Theorem shows  that  $\{T_{nq^{(l)}}\phi_l\}_{n,l \in \mathbb{Z}}$ is a Riesz basis provided  \eqref{s1} holds and $\{T_{nq^{(l)}}\phi_l\}_{n,l \in \mathbb{Z}}$ is linearly independent.

\begin{thm}\label{r3}
	For each $l \in \mathbb{Z}$, let $q^{(l)} \in \mathbb{R}^{+}$. A  nonstationary  frame of translates $\{T_{nq^{(l)}}\phi_l\}_{n,l \in \mathbb{Z}}$ for $L^2(\mathbb{R},dx)$ is a Riesz basis of the space $L^2(\mathbb{R},dx)$ if  $\{T_{nq^{(l)}}\phi_l\}_{n,l \in \mathbb{Z}}$ is linearly independent and \eqref{s1} holds.
\end{thm}
\proof
Assume that $\{T_{nq^{(l)}}\phi_l\}_{n,l \in \mathbb{Z}}$ is linearly independent and (\ref{s1}) holds. Let $s,t \in \mathbb{N}$. By linear independence of $\{T_{nq^{(l)}}\phi_l\}_{n,l \in \mathbb{Z}}$, we can say that $\{T_{nq^{(l)}}\phi_l\}_{n\in \mathbb{F}_1, l\in \mathbb{F}_2}$ is (Riesz) basis for $\mathcal{V}=\text{span}\{T_{nq^{(l)}}\phi_l\}_{n\in \mathbb{F}_1, l\in \mathbb{F}_2}$.  Let $\mathcal{S}_{s,t}$ be its associated frame operator. Then,  dual basis of $\{T_{nq^{(l)}}\phi_l\}_{n\in \mathbb{F}_1, l\in \mathbb{F}_2}$ is of the form $\{\mathcal{S}_{s,t}^{-1}T_{nq^{(l)}}\phi_l\}_{n\in \mathbb{F}_1, l\in \mathbb{F}_2}$. Thus,  every $f \in \mathcal{V}$ can be expressed (uniquely) as
\begin{align*}
	f=\sum_{n\in \mathbb{F}_1}\sum_{l\in \mathbb{F}_2}\langle f , \mathcal{S}_{s,t}^{-1}T_{nq^{(l)}}\phi_l \rangle T_{nq^{(l)}}\phi_l,
	\intertext{and}
	\langle T_{nq^{(l)}}\phi_l , \mathcal{S}_{s,t}^{-1}T_{iq^{(j)}}\phi_j \rangle =\delta_{ni}\delta_{lj}=\begin{cases}
		1, \quad \text{if}\ n=i, \  l=j;\\
		0, \quad \text{otherwise},
	\end{cases}
\end{align*}
for $n,i\in \mathbb{F}_1$ and $l,j\in \mathbb{F}_2$.
Letting $s,t\to \infty$, and using (\ref{s1}), for $n$, $l$, $i$, $j\in \mathbb{Z}$,  we obtain
\begin{align*}
	\langle T_{nq^{(l)}}\phi_l , \mathcal{S}^{-1}T_{iq^{(j)}}\phi_j \rangle =\delta_{ni}\delta_{lj}=\begin{cases}
		1, \quad \text{if}\ n=i, \  l=j; \\
		0, \quad \text{otherwise}.
	\end{cases}
\end{align*}
Thus, $\{T_{nq^{(l)}}\phi_l\}_{n,l \in \mathbb{Z}}$ has a biorthogonal sequence $\{\mathcal{S}^{-1}T_{iq^{(j)}}\phi_j\}_{i,j\in \mathbb{Z}}$. Hence, by Theorem \ref{Riesb},  the sequence $\{T_{nq^{(l)}}\phi_l\}_{n,l \in \mathbb{Z}}$ is a Riesz basis for $L^2(\mathbb{R},dx)$.
\endproof

\section{Linear Independence of  Nonstationary Sequences of Translates}\label{Sec7}
In this section, we discuss linear independence of a nonstationary sequence of translates $\{T_{nq^{(l)}}\phi_l\}_{n,l \in \mathbb{Z}}$ in $L^2(\mathbb{R},dx)$, where   $q^{(l)} \in \mathbb{R}^{+}$  for each $l \in \mathbb{Z}$.  Christensen and Hasannasab in \cite{olehasan} proved an equivalent criteria for a sequence $\{f_k\}_{k \in \mathbb{N}}$ in a separable Hilbert space to be linearly independent in terms of a linear operator on span$\{f_k\}_{k \in \mathbb{N}}$  such that iterated action of that operator on an element $f_1$ in $\mathcal{H}$ gives the sequence $\{f_k\}_{k \in \mathbb{N}}$. They proved the following result.
 \begin{prop}\cite[Proposition 2.3]{olehasan}\label{ooindep}
Consider any sequence $\{f_k\}_{k\in \mathbb{N}}$ in $\mathcal{H}$ for which $span\{f_k\}_{k\in \mathbb{N}}$ is infinite dimensional. Then, the following are equivalent:
\begin{enumerate}[$(i)$]
\item $\{f_k\}_{k\in \mathbb{N}}$ is linearly independent.
\item There exists a linear operator $L\colon \text{span}\{f_k\}_{k\in \mathbb{N}} \to \mathcal{H}$ such that $\{f_k\}_{k=1}^{\infty}=\{L^nf_1\}_{n=0}^{\infty}$.
\end{enumerate}
\end{prop}
Now we give an equivalent criteria of a nonstationary sequence $\{T_{nq^{(l)}}\phi_l\}_{n,l \in \mathbb{Z}}$ to be linearly independent in terms of a linear operator defined on span$\{T_{nq^{(l)}}\phi_l\}_{n,l \in \mathbb{Z}}$ such that iterative action of that operator on $T_{nq^{(1)}}\phi_1$ gives the whole sequence $T_{nq^{(l)}}\phi_l$ when $\text{span}\{T_{nq^{(l)}}\phi_l\}_{n,l \in \mathbb{N}}$ is infinite dimensional. This extends Proposition \ref{ooindep} to  nonstationary sequences of translates  in the space $L^2(\mathbb{R},dx)$.
\begin{thm}\label{indep}
Suppose $\{T_{nq^{(l)}}\phi_l\}_{n,l \in \mathbb{Z}}$ is a  nonstationary sequence of translates  in $L^2(\mathbb{R},dx)$ such that $\text{span}\{T_{nq^{(l)}}\phi_l\}_{n,l \in \mathbb{N}}$ is infinite dimensional. Then, the following conditions are equivalent.
\begin{enumerate}[$(i)$]
\item  $\{T_{nq^{(l)}}\phi_l\}_{n,l \in \mathbb{Z}}$ is linearly independent.\label{eqi0}
\item There exists a linear operator $L \colon \text{span}\{T_{nq^{(l)}}\phi_l\}_{n,l \in \mathbb{Z}} \rightarrow L^2(\mathbb{R},dx)$ such that
\begin{align*}
\big\{L^{m,k}T_{q^{(1)}}\phi_1\big\}_{m,k\in \mathbb{Z}}= \big\{T_{(m+1)q^{(k+1)}}\phi_{k+1}\big\}_{m,k \in \mathbb{Z}}.
\end{align*}\label{eqii}
\end{enumerate}
\end{thm}

\proof
\textcolor[rgb]{1.00,0.00,0.00}{For convenience, write} $\phi_{n,l}=T_{nq^{(l)}}\phi_l$, \  $n$, $l\in \mathbb{N}$.

 $\eqref{eqi0}\implies \eqref{eqii}$: Assume that $\{\phi_{n,l}\}_{n,l\in \mathbb{Z}}$ is linearly independent. Hence, every subset of $\{\phi_{n,l}\}_{n,l\in \mathbb{Z}}$ is also linearly independent. For a fixed $l \in \mathbb{Z}$, define
\begin{align*}
L^{1,0}\phi_{n,l}=
\phi_{n+1,l}, \ n  \in \mathbb{Z}.
\end{align*}
Extend $L^{1,0}$ to an operator on span$\{\phi_{n,l}\}_{n\in \mathbb{N}}$ linearly. Then, for a fixed $l \in \mathbb{Z}$, we have
\begin{align}\label{i1}
\{L^{m,0}\phi_{1,l}\}_{m \in \mathbb{Z}}=
\{\phi_{m+1,l}\}_{m \in \mathbb{Z}}.
\end{align}
Similarly, for fixed $n \in \mathbb{Z}$, we have
\begin{align}\label{i2}
\{L^{0,k}\phi_{n,1}\}_{k \in \mathbb{Z}}=
\{\phi_{n,k+1}\}_{k \in \mathbb{Z}}.
\end{align}
Using $(\ref{i1})$ and $(\ref{i2})$, there exists a linear operator on span$\{\phi_{n,l}\}_{n,l\in \mathbb{N}}$ such that
\begin{align*}
\{L^{m,k}\phi_{1,1}\}_{m,k\in \mathbb{Z}}=
\{\phi_{m+1,k+1}\}_{m,k \in \mathbb{Z}}.
\end{align*}
Thus, $\eqref{eqi0}$ implies $\eqref{eqii}$.

 $\eqref{eqii} \implies \eqref{eqi0}$:  By \eqref{eqii}  there exists a linear operator $L \colon \text{span}\{\phi_{n,l}\}_{n,l \in \mathbb{N}} \rightarrow L^2(\mathbb{R},dx)$ such that
\begin{align*}
\{L^{m,k}\phi_{1,1}\}_{m,k\in \mathbb{Z}}=
\{\phi_{m+1,k+1}\}_{m,k \in \mathbb{Z}}.
\end{align*}

Assume on the contrary that $\{\phi_{n,l}\}_{n,l \in \mathbb{Z}}$ is linearly dependent. Then,  there exist  $R$, $S>0$ such that
\begin{align*}
\sum_{n=-R}^{R}\sum_{l=-S}^{S}c_{n,l}\phi_{n,l}=0
\end{align*}
for some coefficients  $\{c_{n,l}\}_{n \in \{-R,\dots,R\},l\in \{-S,\dots,S\}}$.
Without loss of generality, assume that $c_{R,S}\neq 0$. Then, we can write
\begin{align}
\phi_{R,S}= \sum_{n=-R}^{R-1}\sum_{l=-S}^{S-1}c'_{n,l}\phi_{n,l} + \sum_{n=-R}^{R-1}c'_{n,S}\phi_{n,S} + \sum_{l=-S}^{S-1}c'_{R,l}\phi_{R,l}
\end{align}
for some coefficients $c'_{n,l}$.

Write
\begin{align*}
V:=\text{span}\Big[\{\phi_{n,l}\}_{n \in \{-R,\dots,R-1\}, \atop l\in \{-S,\dots,S-1\}} \bigcup \{\phi_{n,S}\}_{n \in \{-R,\dots,R-1\}} \bigcup \{\phi_{R,l}\}_{l \in \{-S,\dots,S-1\}}\Big].
 \end{align*}
Then, $\phi_{R,S}\in V$. Now, for any $v\in V$, we have
\begin{align*}
L^{1,0}v= L^{1,0}\bigg(\sum_{n=-R}^{R-1}\sum_{l=-S}^{S-1}d_{n,l}\phi_{n,l} + \sum_{n=-R}^{R-1}d'_{n,S}\phi_{n,S} + \sum_{l=-S}^{S-1}d''_{R,l}\phi_{R,l}\bigg)
\end{align*}
for some coefficients $d_{n,l},d'_{n,S}$ and $d''_{R,l}$. Therefore, for any $v\in V$,  we have
\begin{align*}
L^{1,0}v=
\sum\limits_{n=-R+1}^{R}\sum\limits_{l=-S}^{S-1}d_{n,l}\phi_{n,l} + \sum\limits_{n=-R+1}^{R}d'_{n,S}\phi_{n,S} + \sum\limits_{l=-S}^{S-1}d''_{R,l}\phi_{R,l},
\end{align*}
Thus, $L^{1,0}v\in V$.  This implies that $V$ is invariant under $L^{1,0}$. Therefore, $L^{m,0}v\in V$ for $m \in \mathbb{N}\cup \{0\}$. Similarly, $L^{0,k}v\in V$ for $k \in \mathbb{N} \cup \{0\}$. Now for any  $v \in V$ and  for any $m,k \in \mathbb{N}\cup \{0\}$,  we have
\begin{align*}
L^{m,k}v &= L^{m,0}(L^{0,k}v) \\
&= L^{m,0}(v') \in V \ \  \big(\text{where} \ v' = L^{0,k}v \in V\big).
\end{align*}
Thus, \eqref{eqii} implies that $V=\text{span}\{\phi_{n,l}\}_{n \in \{-R,\dots, \infty \}, \atop l\in \{-S,\dots, \infty\}}$, which is a contradiction. Hence, $\{\phi_{n,l}\}_{n,l\in \mathbb{Z}}$ is linearly independent.
This completes the  proof.
\endproof

\begin{rem}
Proposition \ref{ooindep} is a particular case of Theorem \ref{indep}.
Indeed, for a given $n\in \mathbb{N}$ and for $q^{(l)}= c = \text{constant}$, $\l \in \mathbb{N}$, take $nq^{(l)}=a$, $k=0$ and $m \in \mathbb{N} \cup \{0\}$. Then,   $\{f_l\}_{l \in \mathbb{N}}:=\{T_a\phi_l\}_{l\in \mathbb{N}}$ in Proposition \ref{ooindep}.
\end{rem}

$$\textbf{\text{Declaration of competing interest}}$$
There is no competing interest.

\end{document}